\documentclass[12pt]{article}
\usepackage{amsmath} \usepackage{amsfonts} 
\usepackage{natbib}

\begin{document}

\makeatletter	   
\renewcommand{\ps@plain}{%
     \renewcommand{\@oddhead}{\textrm{}\hfil\textrm{\thepage}}%
     \renewcommand{\@evenhead}{\@oddhead}%
     \renewcommand{\@oddfoot}{}
     \renewcommand{\@evenfoot}{\@oddfoot}}
\makeatother     

\newtheorem{theorem}{Theorem}
\newtheorem{lemma}[theorem]{Lemma}
\newtheorem{proposition}[theorem]{Proposition}
\newtheorem{corollary}[theorem]{Corollary}

\title{On rates of convergence for posterior distributions under misspecification }         
\author{By HENG LIAN}        
\date{\small{182 George St, Division of Applied Mathematics, Brown University, \\Providence, RI, 02912 USA.\\Heng\_Lian@brown.edu}}          
\maketitle

\pagestyle{plain}
\begin{center}
\textbf{SUMMARY }
\end{center}

 We extend the approach of Walker (2003, 2004) to the case of misspecified models. A sufficient condition for establishing rates of convergence is given based on a key identity involving martingales, which does not require construction of tests. We also show roughly that the result obtained by using tests can also be obtained by our approach, which demonstrates the potential wider applicability of this method.

\vspace{0.1in}
\noindent\textit{Some key words}: $\alpha$-covering; Bayesian nonparametrics; Prior misspecification. 
\section{INTRODUCTION}       
Bayesian inference distinguishes itself from the frequentist school by its explicit quantification of uncertainty of the parameter with prior specification. The classical approach uses subjective priors elicited from field experts and domain knowledge. The modern Bayesian school have instead shifted attention more towards the construction of priors using formal rules in the hope of dealing with arbitrariness of the prior. 

Asymptotics for infinite-dimensional Bayesian statistics has been receiving a lot of attention recently. In these studies, the Bayesian inference is approached from a frequentist point of view, that is, we assume there is a true underlying probability distribution that generates the data. Naturally one desired property is that as more and more observations are made from the underlying generating mechanism, we will obtain accurate estimate of the true distribution. While traditional Bayesians do not believe in such an assumption, it is shown by Blackwell \& Dubins (1962) that this property is the same as intersubjective agreement, which means two Bayesians will eventually come to roughly the same conclusion after seeing enough data.

The posterior distribution typically behaves well under regular parametric models. Doob showed that consistency is achieved under almost no assumptions on the model, except for a zero measure set under the prior, although in topological terms this set can be large. For infinite-dimensional models, however, the matter is more subtle. Strange behavior can be observed under some priors as documented in Diaconis \& Freeman (1986). 
Given the prior $\Pi$ on the set $\mathcal{P}$ of probability distribution, the posterior is a random measure:
\[ 
\Pi^n(B|X_1,\ldots,X_n)=\frac{\int_B\prod_{i=1}^np(X_i)d\Pi(P)}{\int\prod_{i=1}^np(X_i)d\Pi(P)}
\]
For ease of notation, we will omit the conditioning and only write $\Pi^n(B)$ for the posterior distribution. We say that the posterior is consistent if 
\[
\Pi^n(P\in \mathcal{P}: d(P,P_0)>\epsilon)\rightarrow 0 \mbox{ in $P_0^n$ probability}.
\]
where $P_0$ is the true distribution and $d$ is some suitable distance function between probability measures.

To study rates of convergence, let $\epsilon_n$ be a sequence decreasing to zero, we say the rate is at least $\epsilon_n$ if for sufficiently large constant $M$
\[ 
\Pi^n(P: d(P,P_0)\ge M\epsilon_n)\rightarrow 0 \mbox{ in $P_0^n$ probability}.
\]
We can also have a slightly weaker definition of rates of convergence by replacing $M$ with a sequence $M_n$ and requiring that the above posterior mass converge to zero for any sequence $M_n$ that diverges to infinity.

On the positive side, Schwartz (1965) shows consistency for specific distributions by constructing a sequence of tests of the true distribution against distributions some positive distance away. The tests can trivially be constructed for weak neighborhoods. The construction of similar tests for stronger topology (typically measured in Hellinger distance, for example) is not so straightforward and requires extra works. Barron et al. (1999) gives sufficient conditions that guarantee consistency of infinite-dimensional models by bounding the likelihood ratio under bracketing entropy constraint on sieves. Shen \& Wasserman (2001) studied rates of convergence. A related approach by constructing a sequence of tests appeared in Ghosal et al. (2000). 

The conditions imposed in the above are sufficient but not necessary. It is important to see to what extent these conditions can be relaxed. Another line of work parallel to the development above by Stephen Walker and his collaborators proves consistency and rates of convergence under slightly less stringent conditions. These results are established by constructing a certain supermartingale and consistency and rates of convergence is shown by focusing on the distance of certain predictive distributions to the true one. This approach does not require construction of sequence of tests or sieves. It is shown that this new approach can lead to somewhat weaker sufficient conditions or faster rates.

In Kleijn \& van der Vaart (2006), the authors consider the situation where one cannot expect to achieve consistency since the prior is misspecified. In this case, it is not surprising that the posterior will converge to the distribution in the support of the prior that is closest to the true distribution measured in Kullback-Leibler divergence. Instead of using the usual entropy number or its local version, they used a new concept called covering number for testing under misspecification and studied rates by constructing a sequence of tests between the true distribution $P_0$ and another measure that is not necessarily a probability distribution. The new entropy number can be reduced to the usual entropy in the well-specified case. In this paper, we study the posterior distribution also under the misspecified situation, without constructing a sequence of tests. 

The goal of this paper is two fold. First, we show that the approach in Walker(2003, 2004) can be extended to the situation of misspecified prior rather straightforwardly, by introducing an $\alpha$-entropy condition that is slightly stronger than that of Kleijn \& van der Vaart (2006). Second, we show that using a more refined analysis, a result similar to Theorem 2.2 in Kleijn \& van der Vaart (2006) can be recovered. In particular, it shows that under the well-specified case, this approach indeed is more general than the approach of constructing a  sequence of tests. 

In \S 2, we introduce necessary notations and concepts and present the martingale construction due to Walker (2003). In \S 3, we prove the main result and show that this approach is somehow more general than the one presented in Kleijn \& van der Vaart (2006). We end this paper with a discussion in \S 4.

\section{PRELIMINARIES}

Let $\{X_1, X_2, \ldots\}$ be independent samples generated from distribution $P_0$, with corresponding lower case letter $p_0$ denoting the density with respect to some dominating measure $\mu$. We are given a  collection of distributions $\mathcal{P}$, and a prior $\Pi$ on it with $\Pi(\mathcal{P})=1$. For simplicity, we assume that there exists a unique distribution $P^*\in\mathcal{P}$ that achieves minimum value of Kullback-Leibler divergence to the true distribution, that is 
\begin{eqnarray*}
E_0(\log\frac{p_0}{p^*})\le E_0(\log\frac{p_0}{p}), \; \mbox{for all } p\in\mathcal{P}
\end{eqnarray*} where $E_0$ denotes the expectation under the true distribution $P_0$.

Let $R_n(p)=\prod_{i=1}^n p(X_i)/p^*(X_i)$, then the posterior mass for a set $B$ is 
\begin{eqnarray}\label{ratio}
\Pi^n(B)=\frac{\int_B R_n(p)\Pi(P)}{\int R_n(p)\Pi(P)}
\end{eqnarray}

Following Kleijn \& van der Vaart (2006), for $\epsilon>0, 0<\alpha<1 $ and some suitable semi-metric $d$ on $\mathcal{P}$, we define the $\alpha$-covering of the set $A=\{P\in\mathcal{P}: d(P,P^*)\ge\epsilon\}$ as a collection of convex sets $\{A_1, A_2,\ldots\}$ that covers $A$ with the additional property that for any $j$,
\begin{eqnarray}\label{alphacovering}
\inf_{P\in A_j}-\log E_0(\frac{p}{p^*})^\alpha\ge\frac{\epsilon^2}{4}
\end{eqnarray}
and denote by $N_t(\epsilon,\alpha,A)$ the minimum integer $N$ such that there exists $\{A_1, \ldots, A_N\}$ that forms such a cover, if $N$ is finite.

This condition appears to be stronger than the concept of covering for testing under misspecification introduced by Kleijn \& van der Vaart (2006), which only requires that 
\begin{eqnarray}\label{covering}
\inf_{P\in A_j}\sup_{0<\alpha<1}-\log E_0(\frac{p}{p^*})^\alpha\ge\frac{\epsilon^2}{4}
\end{eqnarray}
In all the examples they gave in their paper, though, we can find a certain value of $\alpha$ only depending on the specification of the model that satisfies our condition. As shown in Kleijn \& van der Vaart (2006), when $\mathcal{P}$ is convex, we have $d^2(P,P^*)\le -\log E_0(p/p^*)^{1/2}$ where $d$ is a generalized Hellinger distance defined by $d^2(P_1,P_2)=\frac{1}{2}\int(p_1^{1/2}-p_2^{1/2})^2p_0/p^*\,d\mu$, which reduces to the usual Hellinger distance in the well-specified case. In this situation, the $1/2$-covering for testing can be replaced by the usual covering as shown in Kleijn \& van der Vaart (2006). In general, allowing $\alpha$ to be different than $1/2$ is required, since in the misspecified case, we cannot guarantee that $-\log E_0(p/p^*)^{1/2}>0$, and we are obliged to choose some smaller $\alpha$ in order to find the covering.

The predictive density constrained to a general set $A$ is defined as 
\begin{eqnarray*}
p_{nA}(x)=\int_A p(x)\Pi^n_A(P)
\end{eqnarray*}, where $\Pi_A^n(P)=1_{\{P\in A\}}\Pi^n(P)/\Pi^n(A)$ is the posterior measure conditioned on $A$. The key identity noted by Walker (2003) is the following:
\begin{eqnarray*}
\int_A R_{n+1}(p)\Pi(P)=\frac{p_{nA}(X_{n+1})}{p^*(X_{n+1})}\int_A R_{n}(p)\Pi(P)
\end{eqnarray*}
as can be verified easily. This in turn implies that 
\begin{eqnarray}\label{key}
E_0[(\int_A R_{n+1}(p)\Pi(P))^\alpha|X_1,\ldots,X_n]=(\int_A R_{n}(p)\Pi(P))^\alpha E_0(\frac{p_{nA}}{p^*})^\alpha
\end{eqnarray}
which means that $\int_A R_{n}(p)\Pi(P)$ is a supermartingale when $E_0(p_{nA}/p^*)^\alpha<1$.

\section{RATES OF CONVERGENCE}
To study rates of convergence, for a sequence $\epsilon_n\rightarrow 0$, we let $A_{n}=\{P\in\mathcal{P}: d(P,P^*)\ge M\epsilon_n\}$, and let $A_{n,j}$ be an $\alpha$-covering of $A_n$, i.e., $\{A_{n,j}\}$ are convex sets that covers $A_n$ and
\begin{eqnarray*}
\inf_{P\in A_{n,j}}-\log E_0(\frac{p}{p^*})^\alpha\ge \frac{M^2\epsilon_n^2}{4}
\end{eqnarray*}
Define 
\begin{equation}\label{L}
L_{k,j}^{(n)}=\int_{A_{n,j}} R_k(p)\Pi(P) 
\end{equation}
and 
\begin{equation*}
I_n=\int_{\mathcal{P}}R_n(p)\Pi(P)
\end{equation*}
To obtain a lower bound for $I_n$, which is the denominator in (\ref{ratio}), we also need a condition on the prior mass for a Kullback-Leibler neighborhood of $p^*$, which is defined as
\begin{eqnarray*}
 B(\epsilon, P^*; P_0)=\{P\in\mathcal{P}: -E_0(\log\frac{p}{p^*})\le\epsilon^2, E_0(\log\frac{p}{p^*})^2\le\epsilon^2\}
\end{eqnarray*}

\begin{theorem}\label{th:main}
Assume that $P^*$ is the unique minimizer in $\mathcal{P}$ of the Kullback-Leibler divergence to the true distribution with $E_0(\log (p_0/p^*))<\infty$. For a sequence $\epsilon_n$ such that $\epsilon_n\rightarrow 0$ and $n\epsilon_n^2\rightarrow\infty$, and $A_n, A_{n,j}$ defined as above. If the following conditions hold

1) $e^{-n\epsilon_nK}\sum_j(A_{n,j})^\alpha$ for a sufficiently large constant $K$

2) $\Pi(B(\epsilon_n, P^*; P_0))\ge e^{-Ln\epsilon_n^2}$ for a sufficiently large constant $L$

\noindent then $\Pi^n(P:d(P,P^*)\ge M\epsilon_n)\rightarrow 0$ in $P_0^n$ probability.
\end{theorem}
\textit{Proof}. First we observe that $P_{nA_{n,j}}\in A_{n,j}$ by the convexity of $A_{n,j}$. From the definition of $\alpha$-covering, $\inf_{P\in A_{n,j}}-\log E_0(p/p^*)^\alpha\ge M^2\epsilon_n^2/4$, so the predictive density satisfies
\[\ E_0(\frac{p_{nA_{n,j}}}{p^*})^\alpha\le e^{-M^2\epsilon_n^2/4}\]
Taking expectations in (\ref{key}), with $A$ replaced by $A_{n,j}$, we get 
\begin{equation*}
E_0(L_{k+1,j}^{(n)})^\alpha\le E_0(L_{k,j}^{(n)})^\alpha e^{-M^2\epsilon_n^2/4}
\end{equation*}
and hence 
\begin{eqnarray*}
E_0(L_{n,j}^{(n)})^\alpha\le e^{-nM^2\epsilon_n^2/4}(\Pi(A_{n,j}))^\alpha
\end{eqnarray*}
The posterior distribution can be bounded as follows:
\begin{eqnarray*}
\Pi^n(A_{n})&\le&\sum_j\Pi^n(A_{n,j})\\
&\le&\sum_j [\Pi^n(A_{n,j})]^\alpha=\sum_j\frac{(L_{n,j}^{(n)})^\alpha}{I_n^\alpha}\\
\end{eqnarray*}

Lemma 7.1 in Kleijn \& van der Vaart (2006) shows that when $n\epsilon_n^2\rightarrow\infty$, for every $C>0$, on a set $\Omega_n$ with probability converging to 1, we have $I_n\ge\Pi(B(\epsilon_n,P^*;P_0))e^{-n\epsilon_n^2(1+C)}$, so we can write
\begin{eqnarray*}
E_0(\Pi^n(A_{n}))&=& E_0(\Pi^n(A_{n})1_{\Omega_n})+E_0(\Pi^n(A_{n})1_{\Omega_n^c})\\
&\le& \frac{E_0\sum_j(L_{n,j}^{(n)})^\alpha}{\Pi(B(\epsilon_n,P^*;P_0))^\alpha e^{-\alpha n\epsilon_n^2(1+C)}}+P_0(\Omega_n^c)\\
&\le& e^{-nM^2\epsilon_n^2/4+\alpha n\epsilon_n^2(1+C)+\alpha n\epsilon_n^2L}\sum_j\Pi(A_{n,j})^\alpha+P_0(\Omega_n^c)
\end{eqnarray*}
which converges to zero by condition $1)$ if $M$ is sufficiently large. $\Box$

For a compact set of models $\mathcal{P}$, we can use the trivial bound 
\[\sum_j\Pi(A_{n,j})^\alpha\le N_t(\epsilon_n,\alpha,A_n)\], which gives the following result similar to Theorem 2.1 in  Kleijn \& van der Vaart (2006), while they used a local version of the entropy instead.
\begin{theorem}
If instead of condition 1) in Theorem \ref{th:main}, we assume $N_t(\epsilon_n,\alpha,A_n)\le e^{n\epsilon_n^2}$, then for sufficiently large constant $M$,
\[
\Pi^n(P:d(P,P^*)\ge M\epsilon_n)\rightarrow 0 \mbox{ in probability.}
\]
\end{theorem}

In order to get optimal rate for parametric models, Kleijn \& van der Vaart (2006) used a more refined assumption. In place of condition $2)$ in Theorem \ref{th:main} above, they assumed 
\begin{eqnarray}\label{local}
\frac{\Pi(P: J\epsilon_n<d(P,P^*)<2J\epsilon_n^2)}{\Pi(B(\epsilon_n,P^*;P_0))}\le e^{n\epsilon_n^2J^2/8}
\end{eqnarray} for all natural numbers $n$ and $J$. In order to recover this result, we need a more careful analysis. 

First, we define $A_n^J=\{P\in\mathcal{P}:M_nJ\epsilon_n\le d(P,P^*)<2M_nJ\epsilon_n\}$, with $\alpha-$covering $\{A_{n,j}^J\}$ defined similarly as before with the property:
$\inf_{P\in A_{n,j}^J}-\log E_0(p/p^*)^\alpha\ge M_n^2J^2\epsilon_n^2/4$. Let $\tilde{A}_{n,j}^J=A_{n,j}^J\cap A_n^J$, note that $\tilde{A}_{n,j}^J$ might not be convex even though $A_{n,j}^J$ is constrained to be so. Similarly, we can define $\tilde{L}_{k,j}^{(n),J}$ as in (\ref{L}) with $A_{n,j}$ replaced by $\tilde{A}_{n,j}^J$. It is easy to see that the following still holds:
\[
E_0(\tilde{L}_{k+1,j}^{(n),J})^\alpha=E_0(\tilde{L}_{k,j}^{(n),J})^\alpha E_0(\frac{p_{n\tilde{A}_{n,j}^J}}{p*})^\alpha\le E_0(\tilde{L}_{k,j}^{(n),J})^\alpha e^{-M_n^2J^2\epsilon_n^2/4}
\]
even though $\tilde{A}_{n,j}^J$ might be nonconvex, since $P_{n\tilde{A}_{n,j}^J}$ is still contained in $A_{n,j}^J$ though not necessarily in $\tilde{A}_{n,j}^J$. 

With $A_n^J$ playing the role of $A_n$ before, the same strategy in the proof of Theorem \ref{th:main} can be followed to show that 
\begin{eqnarray*}
E_0(\Pi^n(A_n^J)1_{\Omega_n})&\le& \frac{E_0\sum_j(\tilde{L}_{n,j}^{(n),J})^\alpha}{\Pi(B(\epsilon_n,P^*;P_0))^\alpha e^{-\alpha n\epsilon_n^2(1+C)}} \\
&\le& e^{-nM_n^2J^2\epsilon_n^2/4+\alpha n\epsilon^2(1+C)}\frac{\sum_j\Pi(\tilde{A}_{n,j}^{J})^\alpha}{\Pi (B(\epsilon_n,P^*;P_0))^\alpha}
\end{eqnarray*}

We will use the notation $N^J_t$ to denote the $\alpha$-covering number for $A_n^J$. We are now ready to prove the following:
\begin{theorem}\label{th:cor}
Assume that $P^*$ is the unique minimizer of the Kullback-Leibler divergence to the true distribution with $E_0(\log (p_0/p^*))<\infty$. For a sequence $\epsilon_n$ such that $\epsilon_n\rightarrow 0$ and $n\epsilon_n^2$ bounded away from zero, and $A_n^J, \{A_{n,j}^J\}_{j=1}^{N_t^J}$ defined as above. If the following conditions hold

1) $N_t^J\le e^{n\epsilon_n^2}$ for all $J\ge 1$

2) (\ref{local}) is satisfied

\noindent Then we have
\[
\Pi^n(P:d(P,P^*)\ge M_n\epsilon_n)\rightarrow 0
\]
in probability for any sequence $M_n\rightarrow\infty$
\end{theorem}
\textit{Proof}. We start by writing
\begin{eqnarray}\label{ineq}
E_0(\Pi^n(A_n))&=&\sum_{J=1}^\infty E_0(\Pi^n(A_n^J))\nonumber\\
&\le& \sum_J E_0(\Pi^n(A_n^J)1_{\Omega_n})+P_0(\Omega_n^c)\nonumber\\
&\le& \sum_J e^{-nM_n^2J^2\epsilon_n^2/4+\alpha n\epsilon_n^2(1+C)}\frac{\sum_j\Pi(\tilde{A}_{n,j}^{J})^\alpha}{\Pi(B(\epsilon_n,P^*;P_0))^\alpha}+P_0(\Omega_n^c)
\end{eqnarray}
we can bound the inner sum for each fixed $J$ as
\begin{eqnarray*}
\sum_j\Pi(\tilde{A}_{n,j}^{J})^\alpha\le N_t^J\Pi(A_n^J)^\alpha\le e^{n\epsilon_n^2}\Pi(A_n^J)^\alpha
\end{eqnarray*}
since $\tilde{A}_{n,j}^J\subset A_n^J$ and using condition 1). Plugging this into (\ref{ineq}) and using condition 2), we get
\begin{eqnarray*}
E_0(\Pi^n(A_n))\le\sum_{J\ge 1} e^{-n\epsilon_n^2M_n^2J^2/4+\alpha n\epsilon_n^2(1+C)+n\epsilon_n^2+\alpha n\epsilon_n^2M_n^2J^2/8}+P_0(\Omega_n^c)
\end{eqnarray*}
By Lemma 7.1 of Kleijn \& van der Vaart (2006), $P_0(\Omega_n^c)$ can be made arbitrarily small by choosing $C$ sufficiently large, under the condition that $n\epsilon_n^2$ is bounded away from zero. For any $C$, the sum above converges to zero since $M_n\rightarrow\infty$. $\Box$
\section{DISCUSSION}
We demonstrated that rates of convergence of posterior distribution under misspecification can be established without construction of a sequence of tests. Theorem \ref{th:cor} we derived above is slightly weaker than Theorem 2.2 in Kleijn \& van der Vaart (2006) due to our use of assumption (\ref{alphacovering}), which is stronger than (\ref{covering}). This said, we are not aware of any examples where the weaker condition (\ref{covering}) provides any advantage over (\ref{alphacovering}). In Walker (2007), the authors demonstrated that using the martingale approach can improve on the rates slightly for some problems. Theorem \ref{th:cor} shows that the results by Kleijn \& van der Vaart (2006) is implied by our result, this is precisely true for well-specified problem, while for misspecified problem this is not conclusive due to the reason stated above. Unfortunately, we have not been able to construct an example that this approach provides a faster rate.

The extension to the case that the prior $\Pi$ depends on $n$, and the case that there exists a finite number of points at minimal Kullback-Leibler divergence to the true distribution should be straightforward.


\end{document}